\documentclass{amsart}
\usepackage{graphicx}
\input{diagrams} 

\newtheorem{theorem}{Theorem}[section]
\newtheorem{lemma}[theorem]{Lemma}
\newtheorem{proposition}[theorem]{Proposition}
\newtheorem{corollary}[theorem]{Corollary}

\theoremstyle{definition}
\newtheorem{definition}[theorem]{Definition}

\theoremstyle{remark}
\newtheorem{remark}[theorem]{Remark}

\numberwithin{equation}{section}

\newcommand{\dd}{\tilde{d}}
\newcommand{\U}{{\mathcal{U}}}
\newcommand{\fU}{{\mathcal{V}}}
\newcommand{\D}{\mathcal{D}}

\newcommand {\hide}[1]{}

\def \eop {\hbox{}\nobreak\hfill $\Box$
\par \goodbreak \smallskip}
\newcommand {\junk}[1]{}
 \newtheorem{algorithm}{\sc Algorithm}

\newcommand{\R}{\mathbb{R}}

\newcommand {\C}     {\mbox{\rm C}}

\newcommand{\F}{\mathbb Q}

\newcommand {\Ker}      {\mbox{\rm Ker}}

\newcommand {\Tot} {{\rm Tot}}

\def\addots{\mathinner{\mkern1mu
\raise1pt\vbox{\kern7pt\hbox{.}}
\mkern2mu\raise4pt\hbox{.}\mkern2mu
\raise7pt\hbox{.}\mkern1mu}}

\begin{document}
\title[On Projections of Semi-algebraic Sets]
{
On Projections of Semi-algebraic Sets Defined by Few Quadratic Inequalities
}
\author{Saugata Basu}
\address{School of Mathematics,
Georgia Institute of Technology, Atlanta, GA 30332, U.S.A.}
\email{saugata@math.gatech.edu}
\thanks{The author was supported in part by an NSF Career Award 0133597 and a 
Sloan Foundation Fellowship.}
\author{Thierry Zell}
\address{School of Mathematics,
Georgia Institute of Technology, Atlanta, GA 30332, U.S.A.}
\email{zell@math.gatech.edu}


\begin{abstract}
Let $S \subset \R^{k + m}$ be a compact semi-algebraic set
defined by 
$
P_1 \geq 0, \ldots, P_\ell \geq 0,
$
where 
$
P_i \in \R[X_1,\ldots,X_k,Y_1,\ldots,Y_m],$ and $\deg(P_i) \leq 2, 
1 \leq i \leq \ell.
$
Let $\pi$ denote the standard projection from $\R^{k + m}$ onto $\R^m$.
We prove that for any $q >0$, the sum of the first $q$ Betti numbers
of $\pi(S)$ is bounded by 
$
\displaystyle{
(k + m)^{O(q\ell)}.
}$
We also present an algorithm for computing the
the first $q$ Betti numbers of $\pi(S)$,
whose complexity is
$\displaystyle{
(k+m)^{2^{O(q\ell)}}.
}$ 
For fixed $q$ and $\ell$, both the bounds are polynomial in $k+m$.
\end{abstract}

\maketitle
\section{Introduction}
\label{sec:intro}
Designing efficient algorithms for computing the Betti numbers of 
semi-algebraic sets is one of the outstanding open problems in 
algorithmic semi-algebraic geometry. There has been some recent progress
in this area. It has been  known for a while that the zero-th
Betti number (which is also the number of connected components) of 
semi-algebraic sets can be computed in single exponential time. Very recently,
it has been shown that even the first Betti number, and more generally
the first 
$q$ Betti numbers for any fixed constant $q$, can be computed
in single exponential time~\cite{BPR05,Basu05b}. Since the problem of 
deciding whether a given semi-algebraic set in $\R^k$ is empty or not is
NP-hard, and that of computing its zero-th Betti number is \#P-hard,
the existence of polynomial time algorithms for computing the Betti numbers is
considered unlikely. 

One particularly interesting case is that of semi-algebraic sets defined
by quadratic inequalities. 
The class of semi-algebraic sets defined by quadratic inequalities 
is the first interesting class of semi-algebraic sets
after sets defined by linear inequalities, in 
which case the problem of computing topological information reduces
to linear programming for which (weakly) polynomial time algorithms are 
known. 
From the point of view of computational complexity,
it is easy to see that the Boolean satisfiability problem 
can be posed as the problem of deciding whether a certain semi-algebraic set 
defined by  quadratic inequalities is empty or not. Thus, deciding 
whether such a set is empty is clearly NP-hard and counting its number
of connected components is $\#$P-hard.
However, semi-algebraic sets defined by quadratic inequalities 
are distinguished
from arbitrary semi-algebraic sets in the sense that, 
{\em if the number of  inequalities is fixed}, then the
sum of their Betti numbers is bounded polynomially in the dimension.
The following bound was proved by Barvinok~\cite{Barvinok97}.

\begin{theorem}
\label{the:barvinok}
Let $S \subset \R^k$ be a semi-algebraic set defined by the inequalities,
$P_1 \geq 0,\ldots,P_\ell \geq 0$, $\deg(P_i) \leq 2, 1 \leq i \leq \ell$.
Then,
$\displaystyle{
\sum_{i=0}^k b_i(S) \leq k^{O(\ell)},
}
$
where $b_i(S)$ denotes the $i$-th Betti number, which is the dimension
of the $i$-th singular cohomology group of $S$, $H^i(S;\mathbb{Q})$, 
with coefficients in $\mathbb{Q}$.
\end{theorem}

In view of Theorem~\ref{the:barvinok}, it is natural to consider
the {\em class of semi-algebraic sets defined by  a fixed number
of quadratic inequalities} from a computational point of view.
Algorithms for computing various topological properties of this
class of semi-algebraic sets have been developed,
starting from the work of  Barvinok~\cite{Barvinok93}, 
who described a polynomial time algorithm for testing emptiness of a set 
defined by a constant number of quadratic inequalities. This was
later generalized and made constructive by Grigoriev and Pasechnik
in~\cite{GP}, where an algorithm is described for computing sample 
points in every connected component of a semi-algebraic set 
defined over a quadratic map. More recently, polynomial time
algorithms have been designed for computing 
the Euler-Poincar\'e characteristic~\cite{Basu05c} as well as
all the Betti numbers~\cite{Basu05a}
of sets defined by a fixed number of quadratic inequalities 
(with different dependence on the number of inequalities in the
complexity bound).
Note also that the problem of deciding the emptiness of a set defined
by a {\em single quartic equation} is already NP-hard and hence
it is unlikely that there exists polynomial time algorithms for
any of the above problems if the degree is allowed to be greater than
two.

A case of intermediate complexity between semi-algebraic sets defined by 
polynomials of higher degree and sets defined by a fixed number of
quadratic sign conditions is obtained by considering 
projections of such sets.
The operation of linear projection of semi-algebraic sets plays a very
significant role in algorithmic semi-algebraic geometry.  It is a
consequence of the Tarski-Seidenberg principle (see for example
\cite{BPR03}, page 61) that the image of a semi-algebraic set under a
linear projection is semi-algebraic, and designing efficient
algorithms for computing properties of projections of semi-algebraic
sets (such as its description by a quantifier-free formula) is a
central problem of the area and is a very well-studied topic (see for
example~\cite{R92} or~\cite{BPR03}, Chapter 14).  However, the
complexities of the best algorithms for computing descriptions of
projections of general semi-algebraic sets is singly exponential in
the dimension and do not significantly improve when restricted to the
class of semi-algebraic sets defined by a constant number of quadratic
inequalities.  
Indeed, any semi-algebraic set can be realized as the projection of a
set defined by quadratic inequalities, and it is not known whether
quantifier elimination can be performed efficiently when the number of
quadratic inequalities is kept constant.  However, we show in this
paper that, with a fixed number of inequalities, the projections of
such sets are topologically simpler than projections of general
semi-algebraic sets.
This suggests, from the point of view of designing efficient
(polynomial time) algorithms in semi-algebraic geometry, that
projections of semi-algebraic sets defined by a constant number of
quadratic inequalities is the next natural class of sets to consider,
after sets defined by linear and (constant number of) quadratic
inequalities, and this is what we proceed to do in this paper.

In this paper, we describe a polynomial time algorithm 
(Algorithm~\ref{alg:main}) for 
computing certain Betti numbers (including the zero-th Betti number
which is the number of connected components)
of projections of sets defined by a constant number of  
quadratic inequalities, without having to compute a semi-algebraic 
description of the projection. 
More precisely, 
let $S \subset
\R^{k+m}$ be a compact semi-algebraic set defined by
$P_1 \geq 0, \ldots, P_\ell \geq 0,$ with $P_i \in
\R[X_1,\ldots,X_k,Y_1,\ldots,Y_m],
\deg(P_i) \leq 2, \; 1 \leq i \leq \ell$. 
Let $\pi:\R^{k+m} \rightarrow \R^m$ be the projection onto the
last $m$ coordinates.  In what follows, the number of inequalities,
$\ell$, used in the definition of $S$ will be considered as some fixed
constant.  Since, $\pi(S)$ is not necessarily describable using only
quadratic inequalities, the bound in Theorem~\ref{the:barvinok} does
not hold for $\pi(S)$ and $\pi(S)$ can in principle be quite
complicated.  Using the best known complexity estimates for quantifier
elimination algorithms over the reals (see~\cite{BPR03}), we get
single exponential (in $k$ and $m$) bounds on the degrees and the
number of polynomials necessary to obtain a semi-algebraic description
of $\pi(S)$.  In fact, there is no known algorithm for computing a
semi-algebraic description of $\pi(S)$ in time polynomial in $k$ and
$m$. Nevertheless, we are able to prove 
that for any fixed constant $q > 0$, the sum of the first $q$
Betti numbers of $\pi(S)$ are bounded by a polynomial in $k$ and $m$.
More precisely, we obtain the following complexity bound (see
Section~\ref{sec:proof}). 

\begin{theorem}
\label{the:main1}
Let $S \subset \R^{k+m}$ be a compact
semi-algebraic set defined by 
\[
P_1 \geq  0, \ldots, P_\ell \geq 0, 
P_i \in \R[X_1,\ldots,X_k,Y_1,\ldots,Y_m], 
\deg(P_i) \leq 2, \; 1 \leq i \leq \ell.
\]
Let $\pi:\R^{k+m} \rightarrow \R^m$ be the projection onto the last
$m$ coordinates. For any $q >0,\; 0 \leq q \leq k$,
\[
\sum_{i=0}^q b_i(\pi(S)) \leq (k + m)^{O(q\ell)}.
\]
\end{theorem}

We also consider the problem of computing the Betti numbers of
$\pi(S)$.  Previously, there was no polynomial time algorithm for
computing any non-trivial topological property of projections of sets
defined by few quadratic inequalities.  We describe a polynomial time
algorithm for computing the first few Betti numbers of $\pi(S)$.
The algorithm (Algorithm~\ref{alg:main} in
Section~\ref{sec:main})
computes $b_{0}(\pi(S)), \ldots, b_q(\pi(S)).$ 
The complexity of the algorithm  is
$\displaystyle{
 (k+m)^{2^{O(q\ell)}}.
}$
If the coefficients of the input polynomials 
are integers  of bit-size bounded by
$\tau$, then the bit-size of the integers
appearing in the intermediate computations and the output
are bounded by
$\tau (k+m)^{ 2^{O(q\ell)}}.$
Note that the output of the algorithm includes
$b_0(\pi(S)),$ which is the number
of connected components of $\pi(S)$.  
Alternatively, one could obtain $b_{0}(\pi(S)), \ldots, b_q(\pi(S))$
by computing a semi-algebraic description of $\pi(S)$ using an
efficient quantifier elimination algorithm (such as the one described
in~\cite{BPR95}) and then using the algorithm described
in~\cite{Basu05b} to compute the first few Betti numbers. However, the
complexity of this method would be worse: single exponential in $k$ and
$m$. Thus, our algorithm is able to
compute efficiently non-trivial topological information about the
projection, even though it does not compute a semi-algebraic
description of that projection (it is not even known whether such a
description could be computed in polynomial time).

In order to obtain Algorithm~\ref{alg:main}, we give a new construction of
a certain spectral sequence, namely
the {\em cohomological descent spectral sequence}, 
which plays a crucial role in the design of the algorithm. 
Even though variants of this spectral sequence have been known for
some time \cite{deligne,stdonat,dugger,houston,murray,vassiliev},
to our knowledge this is the first time it has been used
in designing efficient algorithms.
The new construction that makes this possible is formally analogous to 
that of  the Mayer-Vietoris spectral sequence,
which has been used several times recently 
in designing algorithms for computing 
Betti numbers of semi-algebraic sets (see \cite{B03,Basu05a,Basu05b,BPR05}), 
and thus this new construction %
(see Proposition~\ref{prop:MV} below)
might be of independent interest. 

\section{Main Ideas}
\label{sec:outline}
There are two main ingredients behind the results in this paper.
The first is the use of 
cohomological descent, 
a spectral sequence
first introduced by Deligne~\cite{deligne, stdonat} in the context
of sheaf cohomology. 
This descent spectral sequence is used to compute the
cohomology of the target of a continuous surjection (under certain
hypotheses only, the limit of this spectral sequence is not, in
general, the homology of the target).
The first terms of the sequence are 
cohomology groups 
of certain fibered products over the
surjection, and this allows to bound the Betti numbers of the target
space in terms of the Betti numbers of those fibered products. 
This estimate was first used by Gabrielov, Vorobjov and Zell
in~\cite{GVZ} to give estimates on the Betti numbers of projections of
semi-algebraic sets (and more generally, of semi-algebraic sets
defined by arbitrary quantified formulas) without resorting to
quantifier elimination. Another use of this sequence to establish
upper-bounds can be found in~\cite{Hausdorff} which contains effective
estimates for the Betti numbers of semi-algebraic Hausdorff limits.

The most striking feature of this spectral sequence argument is that 
it enables one to deduce properties (for instance, bounds on 
the Betti numbers) of the projection of a set 
without having to explicitly describe the projection. 
For instance, consider a semi-algebraic subset of $R^k$ defined by a
polynomial having a constant number (say $m$) of monomials (often
referred to as a fewnomial). It is known due to classical results of
Khovansky~\cite{Khovansky} (see also
\cite{B99}) that
the Betti numbers of such sets can be bounded in terms of $m$ and $k$
independent of the degree of the polynomial. 
Using the spectral sequence argument mentioned above, it was proved in 
\cite{GVZ} 
that even the Betti numbers of the projection of such a set
can be bounded in terms of the number of monomials, even though 
it is known (see~\cite{Gabrielov}) that
the projection itself might not admit a description in terms of fewnomials.

The construction of the descent spectral sequence given in~\cite{GVZ}
involves  consideration of join spaces and their filtrations and is 
not directly amenable for algorithmic applications.
In Section~\ref{sec:sequence}, we give an alternate construction of a
descent spectral sequence.
When applied to  
surjections between open sets this spectral sequence
converges to the cohomology of the image. The proof of this fact is
formally analogous to the proof of convergence of the spectral
sequence arising from the generalized Mayer-Vietoris sequence.  
This new proof allows us to identify a certain double complex, whose
individual terms corresponds to the chain groups of the fibered
products of the original set. The fibered product (taken a constant
number of times) of a set defined by few quadratic inequalities is
again a set of the same type.

However, since there
is no known algorithm for efficiently triangulating semi-algebraic sets
(even those defined by few quadratic inequalities) we cannot directly use the
spectral sequence to actually compute the Betti numbers of the projections.
In order to do that we need an additional ingredient.
This second main ingredient is a polynomial time algorithm described in
\cite{Basu05a} for computing a complex whose cohomology groups are
isomorphic to those of a given semi-algebraic set defined by a
constant number of quadratic inequalities. Using this algorithm we are
able to construct a certain double complex, whose associated total
complex is quasi-isomorphic to (implying having isomorphic homology
groups) a suitable truncation of the one obtained from the cohomological
descent spectral sequence mentioned above. This complex is of much
smaller size and can be computed in polynomial time and is enough for
computing the first $q$ Betti numbers of the projection in polynomial
time for any fixed constant $q$.

The rest of the paper is organized as follows.  In
Section~\ref{sec:top_prelim} we recall certain basic facts from
algebraic topology including the notions of complexes, and double
complexes of vector spaces and spectral sequences.  We do not prove
any results since all of them are quite classical and we refer the
reader to appropriate references \cite{Bredon,Mcleary,BPR03} for the
proofs.  In Section~\ref{sec:proof} we prove the estimate on the sum
of Betti numbers  (Theorem~\ref{the:main1}) 
of projections of semi-algebraic sets defined by quadratic
inequalities.  
In Section~\ref{sec:sequence}, we give our
new construction of the cohomological descent spectral sequence 
In Section~\ref{sec:algo}, we briefly describe
Algorithm~\ref{alg:complex} which is used to
compute cohomology groups of semi-algebraic sets given by quadratic
inequalities. This algorithm runs in polynomial time when the number
of inequalities is constant.
We only describe the inputs, outputs and the complexity estimates of
the algorithms, referring the reader to~\cite{Basu05a} 
for more details.  Finally, in Section~\ref{sec:main} we describe our
algorithm (Algorithm~\ref{alg:main}) for computing the first few Betti
numbers of projections of semi-algebraic sets defined by quadratic
inequalities.

\section{Topological Preliminaries}
\label{sec:top_prelim}
We first recall some basic facts from algebraic topology,
related to double complexes, and spectral sequences associated to
double complexes as well as to continuous maps between semi-algebraic sets.
We refer the reader to~\cite{Bredon, Mcleary} for detailed proofs.
We also fix our notations for these objects. All the facts that we
need are well known, and we merely give a brief overview.

\subsection{Complex of Vector Spaces}
A  {\em cochain complex} is a sequence $C^\bullet=\{C^i \mid i \in {\mathbb
Z}\}$ of $\F$-vector spaces together with a sequence of
homomorphisms $\delta^i :C^i \rightarrow C^{i+1}$ for which
$\delta^{i}\circ \delta^{i+1} = 0$ for all $p$.

The cohomology groups, $H^i(C^\bullet)$ are defined by,
\[ 
H^i(C^\bullet) = {Z^i(C^\bullet)}/{B^i(C^\bullet)},
\]
where
$B^i(C^\bullet)= {\rm Im}(\delta^{i-1}),$
and 
$Z^i(C^\bullet) = \Ker(\delta^i).$
The cohomology groups, $H^*(C^\bullet),$ are all
$\F$-vector spaces
(finite dimensional if the vector spaces $C^i$ are themselves finite
dimensional). 
We will henceforth omit reference to 
the field of coefficients $\F$ which is fixed throughout the rest of the
paper.

Given two  complexes, $C^\bullet = (C^i,\delta^i)$ and $D^{\bullet}=
(D^{i},\partial^i)$,
a homomorphism of complexes,
$\phi: C^{\bullet} \rightarrow D^{\bullet},$ is a
sequence of linear maps $\phi^i: C^i \rightarrow D^i$ verifying
$\partial^i\circ \phi^i = \phi^{i+1}\circ\delta^i$ for all $i.$

In other words, the following diagram
is commutative for all $i$.
\[
\begin{array}{ccccccc}
\cdots & \longrightarrow & C^i &
\stackrel{\delta^i}{\longrightarrow} & C^{i+1} & \longrightarrow &\cdots
\\ & &
\Big\downarrow\vcenter{\rlap{$\phi^i$}} & &
\Big\downarrow\vcenter{\rlap{$\phi^{i+1}$}} & & \\ \cdots & \longrightarrow
& D^i &
\stackrel{\partial^i}{\longrightarrow} & D^{i+1} & \longrightarrow &\cdots
\end{array}
\]

A homomorphism of complexes,
$\phi: C^{\bullet} \rightarrow D^{\bullet},$ induces homomorphisms,
$\phi^*: H^*(C^{\bullet}) \rightarrow H^*(D^{\bullet}).$
The homomorphism $\phi$ is called a {\em quasi-isomorphism} if the 
homomorphisms $\phi^*$ are isomorphisms.

\subsection{Double Complexes}
\label{double}
A  double complex
is a bi-graded vector space
$$
{ C^{\bullet,\bullet} = \bigoplus_{i,j \in {\mathbb Z}}
C^{i,j}, } 
$$ 
with co-boundary operators 
$d : C^{i,j} \rightarrow C^{i,j+1}$ and 
$\delta: C^{i,j} \rightarrow C^{i+1,j}$ such that
$d^2=\delta^2=d\delta +\delta d = 0$.
We say that $C^{\bullet,\bullet}$ is a first quadrant
double complex if it satisfies
the condition that $C^{i,j} = 0$ 
when $ij<0$.

Given a double complex $C^{\bullet,\bullet}$, we can construct a 
complex $\Tot^{\bullet}(C^{\bullet,\bullet}),$ called the 
{\em associated total complex of $C^{\bullet,\bullet}$} and defined by
$\displaystyle{ \Tot^n(C^{\bullet,\bullet}) = \bigoplus_{i+j=n}
C^{i,j}, } $ 
with differential 
$\displaystyle{
\D^n : \Tot^{n}(C^{\bullet,\bullet}) 
\longrightarrow \Tot^{n+1}(C^{\bullet,\bullet})
}
$ given by $ D^n= d + \delta$.

{\small
\begin{diagram}
&\vdots &              & \vdots&               &\vdots&              & \\
\luLine&\uTo^{d} &  \luLine            & \uTo^{d}&  \luLine             &\uTo^{d}& \luLine              & \\
\rTo^{\delta}& C^{i-1,j+1} &  \rTo^{\delta}&C^{i,j+1}&\rTo^{\delta}&C^{i+1,j+1}&\rTo^{\delta}&\cdots \\
\luLine&\uTo^{d} &  \luLine            & \uTo^{d}&  \luLine             &\uTo^{d}&        \luLine      & \\
\rTo^{\delta}&C^{i-1,j}  & \rTo^{\delta}& C^{i,j} & \rTo^{\delta} & C^{i+1,j}& \rTo^{\delta}&\cdots \\
\luLine&\uTo^{d} &  \luLine            & \uTo^{d}&   \luLine            &\uTo^{d}&            \luLine  & \\
\rTo^{\delta}& C^{i-1,j-1}  & \rTo^{\delta}& C^{i,j-1} & \rTo^{\delta} & C^{i+1,j-1}& \rTo^{\delta}&\cdots \\
\luLine&\uTo^{d} &   \luLine           & \uTo^{d}& \luLine               &\uTo^{d}&           \luLine   & \\
&\vdots &              & \vdots&               &\vdots&              & \\
\end{diagram}
}

\begin{figure}[hbt] 
\begin{center}
\begin{picture}(0,0)%
\includegraphics{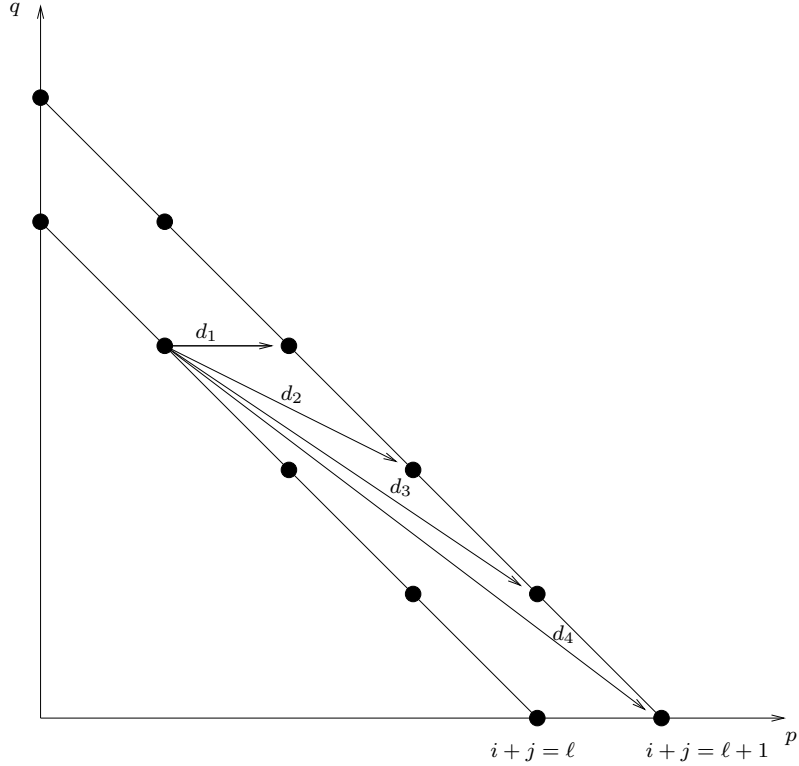}%
\end{picture}%
\setlength{\unitlength}{2565sp}%
\begingroup\makeatletter\ifx\SetFigFont\undefined%
\gdef\SetFigFont#1#2#3#4#5{%
  \reset@font\fontsize{#1}{#2pt}%
  \fontfamily{#3}\fontseries{#4}\fontshape{#5}%
  \selectfont}%
\fi\endgroup%
\begin{picture}(7512,7414)(601,-7394)
\put(6751,-7336){\makebox(0,0)[lb]{\smash{\SetFigFont{8}{9.6}{\familydefault}{\mddefault}{\updefault}{$i + j = \ell+1$}%
}}}
\put(5251,-7336){\makebox(0,0)[lb]{\smash{\SetFigFont{8}{9.6}{\familydefault}{\mddefault}{\updefault}{$i + j = \ell$}%
}}}
\put(8101,-7186){\makebox(0,0)[lb]{\smash{\SetFigFont{8}{9.6}{\familydefault}{\mddefault}{\updefault}{$p$}%
}}}
\put(601,-136){\makebox(0,0)[lb]{\smash{\SetFigFont{8}{9.6}{\familydefault}{\mddefault}{\updefault}{$q$}%
}}}
\put(2401,-3286){\makebox(0,0)[lb]{\smash{\SetFigFont{8}{9.6}{\familydefault}{\mddefault}{\updefault}{$d_1$}%
}}}
\put(3226,-3886){\makebox(0,0)[lb]{\smash{\SetFigFont{8}{9.6}{\familydefault}{\mddefault}{\updefault}{$d_2$}%
}}}
\put(4276,-4786){\makebox(0,0)[lb]{\smash{\SetFigFont{8}{9.6}{\familydefault}{\mddefault}{\updefault}{$d_3$}%
}}}
\put(5851,-6211){\makebox(0,0)[lb]{\smash{\SetFigFont{8}{9.6}{\familydefault}{\mddefault}{\updefault}{$d_4$}%
}}}
\end{picture}

\caption{$d_r: E_r^{i,j} \rightarrow E_r^{i+r, j- r +1}$}
\label{fig:spectral}
\end{center}
\end{figure}

\subsection{Spectral Sequences}
\label{subsec:spectral}
A {\em (cohomology) spectral sequence} is a sequence of bigraded complexes
$\{E_r^{i,j} \mid i,j,r \in \mathbb{Z}, r \geq a\}$
endowed with differentials 
$d_r^{i,j}: E^{i,j}_r \rightarrow E^{i+r,j-r+1}_r$ such that
$(d_r)^2=0$ for all $r.$ 
Moreover, we require the existence of isomorphism between 
the complex $E_{r+1}$ 
and the homology of $E_r$ with respect to $d_r$: 
$$
E_{r+1}^{i,j} \cong H_{d_r}(E_r^{i,j})=
\frac{\ker d_r^{i,j}}{d_r^{i+r,j-r+1}\left(E_r^{i+r,j-r+1}\right)}
$$

The spectral sequence is called a {\em first quadrant spectral sequence}
if the initial complex $E_a$ lies in the first quadrant,
i.e. $E_a^{i,j}=0$ whenever $ij<0.$ In that case, all subsequent
complexes $E_r$ also lie in the first quadrant.  Since the
differential $d_r^{i,j}$ maps outside of the first quadrant for $r>i$,
the homomorphisms of a first quadrant spectral sequence $d_r$ are
eventually zero, and thus the groups $E_r^{i,j}$ are all isomorphic to
a fixed group $E_\infty^{i,j}$ for $r$ large enough, and we say the
spectral sequence is convergent.

Given a double complex $C^{\bullet,\bullet},$ we can associate to it
two spectral sequences, ${'E}_*^{i,j},{''E}_*^{i,j}$ (corresponding to
taking row-wise or column-wise filtrations respectively).

If the double complex lies in the first quadrant, both of these
spectral sequences are first quadrant spectral sequence, and both
converge to $H^*(\Tot^{\bullet}(C^{\bullet,\bullet}))$, meaning
that the limit groups verify 
\begin{equation}\label{eqn:convergence}
\bigoplus_{i+j = n} {'E}_{\infty}^{i,j} \cong 
\bigoplus_{i+j = n} {''E}_{\infty}^{i,j} \cong 
H^n(\Tot^{\bullet}(C^{\bullet,\bullet})),
\end{equation}
for each $n \geq 0$.

The first terms of these are
${'E}_1 = H_{\delta}(C^{\bullet,\bullet}), 
{'E}_2 = H_dH_{\delta}(C^{\bullet,\bullet})$, and
${''E}_1 = H_d (C^{\bullet,\bullet}), 
{''E}_2 = H_\delta H_d (C^{\bullet,\bullet})$.

Given  two (first quadrant) double complexes, $C^{\bullet,\bullet}$ and
$\bar{C}^{\bullet,\bullet},$ a {\em homomorphism of double complexes}
$\phi: C^{\bullet,\bullet} \rightarrow \bar{C}^{\bullet,\bullet}$
is a collection of homomorphisms,
$\phi^{i,j}: C^{i,j} \rightarrow \bar{C}^{i,j},$ such that
the following diagrams commute.

\[
\begin{array}{ccc}
C^{i,j} &
\stackrel{\delta}{\longrightarrow} & C^{i+1,j} \\ 
\Big\downarrow\vcenter{\rlap{$\phi^{i,j}$}} & &
\Big\downarrow\vcenter{\rlap{$\phi^{i+1,j}$}} \\ 
\bar{C}^{i,j} &
\stackrel{\delta}{\longrightarrow} & \bar{C}^{i+1,j}
\end{array}
\]

\[
\begin{array}{ccc}
 C^{i,j} &
\stackrel{d}{\longrightarrow} & C^{i,j+1}\\
\Big\downarrow\vcenter{\rlap{$\phi^{i,j}$}} & &
\Big\downarrow\vcenter{\rlap{$\phi^{i,j+1}$}} \\ 
\bar{C}^{i,j} &
\stackrel{d}{\longrightarrow} & \bar{C}^{i,j+1}
\end{array}
\]

A homomorphism of double complexes, $\phi: C^{\bullet,\bullet}
\rightarrow \bar{C}^{\bullet,\bullet}$ induces homomorphisms
$\phi_r^{i,j}: E_r^{i,j} \rightarrow \bar{E}_r^{i,j}$ between the
terms of the associated spectral sequences (corresponding either to
the row-wise or column-wise filtrations).  

We will need the following useful fact (see~\cite{Mcleary}, page 66,
Theorem 3.4 for a proof).
\begin{theorem}
\label{the:spectral}
If $\phi_s^{i,j}$ is an isomorphism for some $s \geq 1$ (and all
$i,j$), then $E_r^{i,j}$ and $\bar{E}_r^{i,j}$ are isomorphic for all
$r \geq s$. In other words, the induced homomorphism, 
$\phi:\Tot^{\bullet}(C^{\bullet,\bullet}) \longrightarrow
\Tot^{\bullet}(\bar{C}^{\bullet,\bullet})$ is a quasi-isomorphism.
\end{theorem}

\section{Proof of Theorem~\ref{the:main1}}
\label{sec:proof}
The proof of Theorem~\ref{the:main1} relies on the bounds from
Theorem~\ref{the:barvinok}, and on the following theorem that appears
in~\cite{GVZ}. 

\begin{theorem}
\label{the:ss}
Let $X$ and $Y$ be two semi-algebraic sets and $f: X \to Y$ a
semi-algebraic continuous surjection such that $f$ is closed.
Then for any integer $n$, we have

\begin{equation}\label{eqn:ss}
b_n(Y) \leq \sum_{i+j=n} b_j(W^i_f(X)),
\end{equation}
where $W^i_f(X)$ denotes the $(i+1)$-fold fibered product of $X$ over $f$:
\[
W^i_f(X) = \{(\bar{x}_0,\ldots,\bar{x}_i) \in X^{i+1}
\mid f(\bar{x}_0) = \cdots = f(\bar{x}_i)\}.
\]
\end{theorem}

This theorem follows from the existence of a spectral
sequence $E_r^{i,j}$ converging to $H^*(Y)$ and such that $E_1^{i,j}
\cong H^j(W^i_f(X)).$ Since, in any spectral sequence, the dimensions of
the terms $E_r^{i,j}$ are decreasing when $i$ and $j$ are fixed and
$r$ increases, we obtain using the definition~\eqref{eqn:convergence}
of convergence:
\[
b_n(Y)=%
\sum_{i+j=n} \dim \left(E_{\infty}^{i,j}\right)
\leq \sum_{i+j=n} \dim \left(E_1^{i,j}\right),
\]
yielding inequality~\eqref{eqn:ss}.

The spectral sequence $E_r^{i,j}$, known as {\em cohomological
descent}, originated with the work of Deligne~\cite{deligne,
stdonat}, in the framework of sheaf cohomology. In~\cite{GVZ}, the
sequence is obtained as the spectral sequence associated to the
filtration of  
an infinite dimensional topological object, 
the {\em join space}, constructed
from $f$. For the purposes of Algorithm~\ref{alg:main}, we will give a
different construction of this sequence (see Section~\ref{sec:sequence}).

\begin{proof}[Proof of Theorem~\ref{the:main1}:]
Since $S$ is compact, the semi-algebraic continuous surjection $\pi: S
\rightarrow \pi(S)$ is closed: 
applying Theorem~\ref{the:ss} to $\pi$,
inequality~\eqref{eqn:ss} yields
for each $n$ with $0 \leq n \leq q$,
\begin{equation}
\label{eqn:inequality1}
b_n(\pi(S)) \leq \sum_{i+j=n} b_j(W_\pi^i(S)).
\end{equation}
 
Notice that 
$\displaystyle{
W_\pi^i(S) = \{(\bar{x}_0,\ldots,\bar{x}_i,y) \mid P_{h}(\bar{x}_t,y) \geq 0,
1 \leq h \leq \ell, 0 \leq t \leq i \}
}
$.
Thus, each $W_\pi^i(S) \subset \R^{(i+1)k +m}$ is defined by
$\ell(i+1)$ quadratic inequalities.  Applying the bound in
Theorem~\ref{the:barvinok} we get that,
\begin{equation}
\label{eqn:inequality2}
b_j(W_\pi^i(S)) \leq ((i+1)k + m)^{O(\ell(i+1))}.
\end{equation}

Using inequalities~\eqref{eqn:inequality1} and~\eqref{eqn:inequality2}
we get that,
\[
\sum_{i=0}^q b_i(\pi(S)) \leq (k+m)^{O(q\ell)},
\]
which proves the theorem.
\end{proof}

\section{Cohomological Descent}
\label{sec:sequence}

This section is devoted a new construction of the cohomological
descent spectral sequence (already discussed in
Section~\ref{sec:proof}). In Theorem~\ref{the:descent}, we obtain this
sequence as the spectral sequence associated to a double complex
associated to the fibered powers of $X$, rather than through the
filtration of the join space.  Convergence to the cohomology of the
target space occurs when the map $f: X \to Y$ is {\em locally split}
(see definition below). By deformation, we are able to extend the
result to our case of interest: the projection of a compact
semi-algebraic set (Corollary~\ref{cor:descent2}).

We will use this construction for
Algorithm~\ref{alg:main}.

\begin{definition}\label{df:split}
A continuous surjection $f:X \to Y$ is called {\em locally split} if there
exists an open covering $\U$ of $Y$ such that for all $U \in \U$,
there exists a continuous section $\sigma: U \to X$ of $f$,
i.e. $\sigma$ is a continuous map such that $f(\sigma(y))=y$ for all
$y \in U$. 
\end{definition}

In particular, if $X$ is an open semi-algebraic set and $f: X \to Y$
is a projection, the map $f$ is obviously locally split. This specific
case is what we will use in Algorithm~\ref{alg:main}, as we will reduce the
projection of compact semi-algebraic sets to projections of
open semi-algebraic sets (see Proposition~\ref{prop:open}) in
order to apply the spectral sequence.

Recall that for any semi-algebraic surjection $f:X \to Y$, we denoted
by $W_f^p(X)$ the $(p+1)$-fold fibered power of $X$ over $f$,
\[
W^p_f(X) = \{(\bar{x}_0,\ldots,\bar{x}_p) \in X^{p+1}
\mid f(\bar{x}_0) = \cdots = f(\bar{x}_p)\}.
\]
The map $f$ induces for each $p \geq 0,$ 
a map from $W_f^p(X)$ to $Y$,
sending $(\bar{x}_0,\ldots,\bar{x}_p)$ to the common value 
$f(\bar{x}_0) = \cdots = f(\bar{x}_p),$
and 
abusing notations a little
we will denote this  map by $f$ as well. 

\subsection{Singular (co-)homology}
\label{subsec:singular}
We recall here the basic definitions related to singular (co-)homology
theory directing the reader to \cite{Hatcher} for details.
 
For any semi-algebraic set $X$, let
$C_{\bullet}(X)$ denote the complex of singular chains
of $X$ with boundary map denoted by $\partial$.

Recall that $C_{\bullet}(X)$ is defined as follows:
For $m \geq 0$, a singular $m$-simplex $s$ is a continuous map,
$s: \Delta_m \rightarrow X$, where $\Delta_m$ is the
standard $m$-dimensional simplex defined by,
\[
\Delta_m = \{(t_0,\ldots,t_m)\;\mid\; t_i \geq 0, \sum_{i=0}^m t_i = 1 \}.
\]
$C_m(X)$ is the vector space spanned by all singular $m$-simplices
with boundary maps defined as follows. 
As usual we first define the face maps 
\[f_{m,i}: \Delta_{m} \rightarrow \Delta_{m+1},\]
by
$f_{m,i}((t_0,\ldots,t_m)) = (t_0,\ldots,t_{i-1},0,t_{i+1},\ldots,t_{m+1})$. 

For a singular $m$-simplex $s$ we define 
\begin{equation}
\label{eqn:partialdef} 
\partial s = \sum_{i=0}^{m} (-1)^i s\circ f_{m-1,i}.
\end{equation}
and extend $\partial$ to $C_m(X)$ by linearity.
We will denote by $C^\bullet(X)$ the dual complex
and by $d$ the corresponding co-boundary map. More precisely,
given $\phi \in C^m(X)$, and a singular $(m+1)$-simplex $s$ of $X$,
we have
\begin{equation}
\label{eqn:ddef} 
d \phi(s) = \sum_{i=0}^{m+1} (-1)^i \phi(s \circ f_{m,i}).
\end{equation}

If $f : X \rightarrow Y$ is a continuous map, then it naturally induces a 
homomorphism $f_*: {C}_{\bullet}(X) \rightarrow {C}_{\bullet}(Y)$
by defining, for each singular $m$-simplex $s:\Delta_m \rightarrow X$,
$f_*(s) = s \circ f: \Delta_m \rightarrow Y$, which is a
singular $m$-simplex of $Y$.
We will denote by $f^*:{C}^{\bullet}(Y) \rightarrow {C}^{\bullet}(X)$
the dual homomorphism.
More generally, suppose that 
$s = (s_0,\ldots,s_p):\Delta_m \rightarrow W_f^p(X)$ is a singular
$m$-simplex of $W_f^{p}(X).$
Notice that each component, $s_i, 0 \leq i \leq p$ are themselves
singular $m$-simplices of $X$ and that $f_*(s_0) = \ldots = f_*(s_p)$ 
are equal as singular $m$-simplices of $Y.$ We will denote
their common image by $f_*(s)$.

\bigskip

We will require the notion of small simplices subordinate to an open
covering of a topological space (see~\cite{Hatcher}). Assuming that
$f: X \to Y$ is locally split, let $\U$ be an open covering of $Y$ on
which local continuous sections exist.  We denote by $\fU$ the open
covering of $X$ given by the inverse images of elements of $\U$, i.e.
$\fU=\{f^{-1}(U) \mid U \in \U\}$.  We let $C^\U_\bullet(Y)$ be the
subcomplex of $C_\bullet(Y)$ spanned by those singular simplices of
$Y$ whose images are contained in some element of the cover $\U$.
Similarly, we let $C^\fU_\bullet(X)$ be the subcomplex of
$C_\bullet(X)$ spanned by the simplices of $X$ with image in $\fU$,
and more generally, for any integer $p$, $C^\fU_\bullet(W_f^p(X))$
denotes the subcomplex of $C_\bullet(W_f^p(X))$ spanned by simplices
with image contained in $V^{p+1}$ for some $V\in \fU$.  The
corresponding dual cochain complexes will be denoted by
$C_\U^\bullet(Y)$ and $C_\fU^\bullet(W_f^p(X))$ respectively. 
 We will
 henceforth call any singular simplex of $C^\U_\bullet(Y)$ and any
 singular simplex of $C^\fU_\bullet(W_f^p(X))$ {\em admissible}
 simplices.
The inclusion homomorphism,
${\iota}_\bullet:  C^{\U}_{\bullet}(Y) \hookrightarrow C_\bullet(Y)$ induces 
a dual homomorphism,
${\iota}^\bullet:  C^\bullet(Y) \rightarrow C_{\U}^{\bullet}(Y)$. 
We also have corresponding induced homomorhisms,
${\iota}^{\bullet}: C^\bullet(W_f^p(X))\rightarrow 
C_{\fU}^\bullet(W_f^p(X))$ for each $p \geq 0$.

\begin{proposition}
\label{prop:hatcher}
The homomorphism ${\iota}^{\bullet}: C^\bullet(Y) \rightarrow C_\U^{\bullet}(Y)$ 
(resp.  $C^\bullet(W_f^p(X))\rightarrow 
C_{\fU}^\bullet(W_f^p(X))$ for each $p \geq 0$) is a chain homotopy
equivalence.
In particular, we have
$H^*(C_\U^\bullet(Y)) \cong H^*(C^\bullet(Y)) \cong  H^*(Y)$
and 
$H^*(C_\fU^\bullet(W_f^p(X)))
\cong
H^*({C}^\bullet(W_f^p(X))) \cong H^*(W_f^p(X)).$
\end{proposition}
\begin{proof}
  This follows from a similar result for homology, see Proposition
  2.21 in~\cite{Hatcher}.
\end{proof}

\subsection{A long exact sequence}
For each $p \geq 0$, we now define a homomorphism,
\[
\delta^p: C^{\bullet}(W_f^p(X)) \longrightarrow 
C^{\bullet}(W_f^{p+1}(X))
\]
as follows:
for each $i, 0 \leq i \leq p$, define
$\pi_{p,i}: W_f^p(X) \rightarrow W_f^{p-1}(X)$  by,
\[
\pi_{p,i}({x}_0,\ldots,{x}_p) = ({x}_0,\ldots,\widehat{{x}_i},
\ldots,{x}_p)
\] 
($\pi_{p,i}$ drops the $i$-th coordinate).

We will denote by $(\pi_{p,i})_*$ the induced map on
$C_\bullet(W_f^{p}(X)) \rightarrow C_\bullet(W_f^{p-1}(X))$
and let $\pi_{p,i}^*: C^{\bullet}(W_f^{p-1}(X)) \rightarrow 
C^{\bullet}(W_f^{p}(X))$
denote the dual map.
For $\phi \in C^{\bullet}(W_f^p(X))$, we define
$\delta^p \;\phi$ by,
\begin{equation}\label{eqn:deltadef}
\delta^p\;\phi =\sum_{i=0}^{p+1}(-1)^i \pi_{p+1,i}^*\ \phi.
\end{equation}
Note that for any open covering $\fU$ of $X$, the map $\delta^p$
induces by restriction a map
$C^{\bullet}_\fU(W_f^p(X)) \to
C^{\bullet}_\fU(W_f^{p+1}(X))$ which we will still denote by $\delta^p$.

The following proposition is analogous to the exactness of the
generalized Mayer-Vietoris sequence
(cf. Lemma 1 in~\cite{B03}).

\begin{proposition}
\label{prop:MV}
Let $f:X \to Y$ be a continuous, locally split surjection, where $X$
and $Y$ are
semi-algebraic
subsets of $\R^n$ and $\R^m$ respectively. Let $\U$ denote an open
covering of $Y$ in which continuous sections of $f$ can be defined on
every $U \in \U$, and  $\fU$ denote the open covering of $X$
obtained by inverse image of $\U$ under $f$.
The following sequence is exact.
\[
0 \longrightarrow C_\U^\bullet(Y) \stackrel{f^*}{\longrightarrow}
C_\fU^\bullet(W_f^0(X)) \stackrel{\delta^0}{\longrightarrow}
C_\fU^\bullet(W_f^1(X)) \stackrel{\delta^1}{\longrightarrow}
\cdots
\]
\[
\cdots
\stackrel{\delta^{p-1}}{\longrightarrow}
C_\fU^\bullet(W_f^p(X)) \stackrel{\delta^{p}}{\longrightarrow}
C_\fU^\bullet(W_f^{p+1}(X)) \stackrel{\delta^{p+1}}{\longrightarrow}
\cdots
\]
\end{proposition}

\begin{proof}
We will start by treating separately the first two positions in the
sequence, then prove exactness for $p\geq 1$.
\begin{enumerate}
\item $f^*: C_\U^\bullet(Y) \to C_\fU^\bullet(X)$ is injective.

\medskip

Let $U \in \U$ and let $s$ be a simplex whose image is contained in
$U$. If $\sigma$ is a continuous section of $f$ defined on $U$, the
simplex $t=\sigma_*(s)$ is in $C^\fU_\bullet(X)$, and verifies
$f_*(t)=s$. Hence, $f_*: C^\fU_\bullet(X) \to C^\U_\bullet(Y)$ is
surjective, so $f^*$ is injective.

\medskip

\item $f^*(C_\U^\bullet(Y))=\ker \delta^0$.

\medskip

Let $\phi \in C_\fU^m(X)$. Any simplex $s\in C^\fU_m(W^1_f(X))$ is a
pair $(s_0,s_1)$ of simplices in $C^\fU_m(X)$ verifying
$f_*(s_0)=f_*(s_1)$. We then have
$\delta^0\phi(s)=\phi(s_1)-\phi(s_0)$.
If $\phi=f^*\psi$ for some $\psi\in C_\U^m(Y),$ we have for any $s$,
\begin{equation*}
\delta^0\phi(s)=f^*\psi(s_1)-f^*\psi(s_0)=\psi(f_*(s_1))-\psi(f_*(s_0))=0,
\end{equation*}
since we must have $f_*(s_0)=f_*(s_1)$. Thus, we have
$f^*(C_\U^\bullet(Y))\subset\ker \delta^0$.

\medskip

Conversely, if $\phi$ is such that $\delta^0\phi=0$, this means that
for any pair $(s_0,s_1)$ of simplices in $C^\fU_m(X)$ verifying
$f_*(s_0)=f_*(s_1)$, we have $\phi(s_0)=\phi(s_1)$. Since we just
proved in part~(A) that $f_*: C^\fU_\bullet(X) \to C^\U_\bullet(Y)$ is
surjective, any element $t\in C^\U_m(Y)$ is of the form $t=f_*(s)$ for some
$s \in C^\fU_m(X)$. Thus, we can define $\psi \in C_\U^m(Y)$ by
$\psi(t)=\phi(s)$, and the condition $\delta^0\phi=0$ ensures that
$\psi$ is well defined since its value does not depend on the choice
of $s$ in the representation $t=f_*(s)$. This yields the reverse
inclusion, and hence exactness at $p=0$.

\item
$\delta^{p+1}\circ\delta^{p} = 0$: \\
From the definitions of the maps $\pi_{p+1,i}^*, \pi_{p+2,j}^*$ 
we have that for $0 \leq i \leq p+1, 0 \leq j \leq p+2,$
\begin{equation}
\label{eqn:exact}
\pi_{p+2,j}^*\circ\pi_{p+1,i}^*(\phi) = 
\pi_{p+2,i+1}^*\circ\pi_{p+1,j}^*
(\phi)
\;{\rm if}\; j < i. 
\end{equation}

Let ${\phi} \in C_\fU^m(W_f^p(X))$.
Now from the definitions of $\delta^p$ and $\delta^{p+1}$ we have that,

\begin{align*}
\delta^{p+1} \circ \delta^{p}({\phi}) 
&= \delta^{p+1}\left( \sum_{i=0}^{p+1} (-1)^i\pi_{p+1,i}^*(\phi)\right); \\
&=  \sum_{i=0}^{p+1} (-1)^i \delta^{p+1}(\pi_{p+1,i}^*(\phi)); \\
&=  \sum_{i=0}^{p+1} \sum_{j=0}^{p+2} (-1)^{i+j} \pi_{p+2,j}^*
\circ\pi_{p+1,i}^*(\phi); \\
&=  \sum_{i=0}^{p+1} \left[\sum_{0 \leq j < i}(-1)^{i+j} \pi_{p+2,j}^*
\circ\pi_{p+1,i}^*(\phi)
+ \sum_{i \leq j \leq p+2} (-1)^{i+j} 
\pi_{p+2,j}^*\circ\pi_{p+1,i}^*(\phi) \right]; \\
&= \sum_{i \leq  j} (-1)^{i+j} \pi_{p+2,j}^*\circ\pi_{p+1,i}^*(\phi) +
\sum_{i > j}(-1)^{i+j} \pi_{p+2,j}^*\circ\pi_{p+1,i}^*(\phi)).
\end{align*}

Now using Equation~\eqref{eqn:exact}, the previous line becomes
\begin{eqnarray*}%
&=& \sum_{i \leq  j} (-1)^{i+j} \pi_{p+2,j}^*\circ\pi_{p+1,i}^*(\phi) +
\sum_{i > j}(-1)^{i+j} \pi_{p+2,i+1}^*\circ\pi_{p+1,j}^*(\phi)). 
\end{eqnarray*}
Interchanging $i$ and $j$ in the second summand of the previous line,
we get
\begin{eqnarray*}%
&=& \sum_{i \leq  j} (-1)^{i+j} \pi_{p+2,j}^*\circ\pi_{p+1,i}^*(\phi) +
\sum_{i < j}(-1)^{i+j} \pi_{p+2,j+1}^*\circ\pi_{p+1,i}^*(\phi)). 
\end{eqnarray*}
Finally, replacing $j+1$ by $j$ in the second summand above, we obtain
\begin{eqnarray*}%
&=& \sum_{i \leq  j} (-1)^{i+j} \pi_{p+2,j}^*\circ\pi_{p+1,i}^*(\phi) +
\sum_{i < j-1}(-1)^{i+j-1} \pi_{p+2,j}^*\circ\pi_{p+1,i}^*(\phi));  
\end{eqnarray*}
and isolating the terms corresponding to $j=i$ and $j=i+1$ in the
first sum gives
\begin{eqnarray*}
&=& (-1)^{2i}\pi_{p+2,i}^*\circ \pi_{p+1,i}^*(\phi) +
(-1)^{2i+1}\pi_{p+2,i}^*\circ\pi_{p+1,i+1}^*(\phi) \\ 
&&+\sum_{i <  j-1} (-1)^{i+j} \pi_{p+2,j}^*\circ\pi_{p+1,i}^*(\phi) +
\sum_{i < j-1}(-1)^{i+j-1} \pi_{p+2,j}^*\circ\pi_{p+1,i}^*(\phi)); \\
&=& 0;
\end{eqnarray*}
(since, again, by Equation~\eqref{eqn:exact}, we have
$\pi_{p+2,i}^*\circ\pi_{p+1,i+1}^*=\pi_{p+2,i}^*\circ\pi_{p+1,i}^*$).

\item ${\rm Im}(\delta^{p}) \supset {\rm Ker}(\delta^{p+1})$: \\ 
Let $\phi \in {\rm Ker}(\delta^{p+1}).$ 
In  other words, for each admissible 
singular $m$-simplex 
$s= (s_0,\ldots,s_{p+1}):\Delta_m \rightarrow W_f^{p+2}(X)$

\begin{equation}
\label{eqn:exact2}
\sum_{i=0}^{p+2} (-1)^i \phi((s_0,\ldots,\hat{s}_i,\ldots,s_{p+2})) = 0.
\end{equation}

For each admissible 
singular $m$-simplex $s$ of $Y$ 
let $s_*$ denote  
denote a fixed admissible 
singular $m$-simplex of $X$ 
such that $f_*(s_*) = s$. 
Such a choice is possible since, as we proved in part~(A), $f_*$ is
surjective onto $C^\U_\bullet(Y).$ Let $\psi \in C_\fU^m(W_f^p(X))$ be
defined as follows.  For an admissible singular $m$-simplex $t =
(t_0,\ldots,t_p)$ of $W_f^p(X)$ we define
$$
\displaylines{
\psi(t) = \phi(f_*(t),t_0,\ldots,t_p).
}
$$

Now for an admissible 
singular $m$ simplex 
$t = (t_0,\ldots,t_{p+1})$ of $W_f^{p+1}$, with
\begin{align*}
\delta^{p}\psi(t) &= 
\sum_{i=0}^{p+1} (-1)^{i} \pi_{p+1,i}^*\psi(t) \\ 
&= \sum_{i=0}^{p+1} (-1)^i \psi((t_0,\ldots,\hat{t}_i,\ldots,t_{p+1})) \\
&= \sum_{i=0}^{p+1} (-1)^i
\phi((f_*(t),t_0,\ldots,\hat{t}_i,\ldots,t_{p+1})) \\ 
\end{align*}

Now let $s$ denote the admissible 
singular $m$-simplex of 
$W_f^{p+2}(X)$ defined by $s = (f_*(t),t_0,\ldots,\hat{t}_i,\ldots,t_{p+1}).$
Now applying Equation \eqref{eqn:exact2}, we get
$$
\displaylines{
\sum_{i=0}^{p+2} (-1)^i \phi((s_0,\ldots,\widehat{s_i},\ldots,s_{p+2})) = 0.
}
$$

Separating the first term from the rest we obtain,
$$
\displaylines{
\phi((t_0,\ldots,t_{p+1}) = 
\sum_{i=0}^{p+1} (-1)^i \phi((f_*(t),t_0,\ldots,\hat{t}_i,\ldots,t_{p+1})) 
= \delta^{p}\psi(t).
}
$$
\end{enumerate}
This finally proves the exactness of the sequence.
\end{proof}

\subsection{The descent double complex}
Now, let $D^{\bullet,\bullet}(X)$ denote the  double complex defined by,
$D^{p,q}(X) = C^q(W_f^p(X))$ with vertical and horizontal
homomorphisms
given by $\dd^q=(-1)^p d^q$ and $\delta$ respectively,
where $d$ is the singular coboundary operator~\eqref{eqn:ddef} 
and $\delta$ is the map defined in~\eqref{eqn:deltadef}. Also, let
$D^{p,q}(X) = 0$ if $ p < 0$ or $q < 0$. 

\[
\begin{array}{cccccccc}

& & \vdots  && \vdots  && \vdots  & \\
& &
\Big\uparrow\vcenter{\rlap{$\dd$}} & &
\Big\uparrow\vcenter{\rlap{$\dd$}} & &
\Big\uparrow\vcenter{\rlap{$\dd$}} &  \\

0 & \longrightarrow & C^3(W_f^0(X)) &
\stackrel{\delta}{\longrightarrow} & C^3(W_f^1(X)) &
\stackrel{\delta}{\longrightarrow} & 
C^3(W_f^2(X)) &
 \longrightarrow 
\\

& &
\Big\uparrow\vcenter{\rlap{$\dd$}} & &
\Big\uparrow\vcenter{\rlap{$\dd$}} & &
\Big\uparrow\vcenter{\rlap{$\dd$}} &  \\

0 & \longrightarrow & C^2(W_f^0(X)) &
\stackrel{\delta}{\longrightarrow} & 
C^2(W_f^1(X)) &
\stackrel{\delta}{\longrightarrow} & 
C^2(W_f^2(X)) &
 \longrightarrow 
\\
& &
\Big\uparrow\vcenter{\rlap{$\dd$}} & &
\Big\uparrow\vcenter{\rlap{$\dd$}} & &
\Big\uparrow\vcenter{\rlap{$\dd$}} &  \\
0 & \longrightarrow & C^1(W_f^0(X)) &
\stackrel{\delta}{\longrightarrow} & 
C^1(W_f^1(X)) &
\stackrel{\delta}{\longrightarrow} & 
C^1(W_f^2(X)) &
 \longrightarrow 
\\

& &
\Big\uparrow\vcenter{\rlap{$d$}} & &
\Big\uparrow\vcenter{\rlap{$d$}} & &
\Big\uparrow\vcenter{\rlap{$d$}} &  \\

0 & \longrightarrow & C^0(W_f^0(X)) &
\stackrel{\delta}{\longrightarrow} & 
C^0(W_f^1(X)) &
\stackrel{\delta}{\longrightarrow} & 
C^0(W_f^2(X)) &
 \longrightarrow 
\\

& &
\Big\uparrow\vcenter{\rlap{$d$}} & &
\Big\uparrow\vcenter{\rlap{$d$}} & &
\Big\uparrow\vcenter{\rlap{$d$}} &  \\

& & 0 && 0 && 0 & \\
\end{array}
\]

\begin{lemma}
\label{lem:double}
The families of maps $\dd$ and $\delta$ make $D^{\bullet,\bullet}$ into
a double complex.
\end{lemma}

\begin{proof}
We need to check that $\dd^2=\delta^2=\dd\delta+\delta \dd=0$. We know
that $\dd^2=d^2=0$ since $C^\bullet(W^p_f(X))$ is a cochain
complex for all $p$, and we proved that $\delta^2=0$ in
Proposition~\ref{prop:MV}.

Now, suppose that ${\phi} \in C^q(W^p_f(X))$ and let
$s = (s_0,\ldots,s_{p+1})$ 
be an admissible 
singular $(q+1)$-simplex of $W_f^{p+1}(X)$.
Then,
\begin{align*}
\dd(\delta {\phi})(s)
&= \dd\left(\sum_{i=0}^{p+1} (-1)^i {\phi}(
(s_0,\ldots,\widehat{s_i},\ldots,s_{p+1}))\right);\\
&=(-1)^p \sum_{j=0}^{q+1} \sum_{i=0}^{p+1} (-1)^{i+j} 
{\phi}(s_0 \circ f_{q,j},\ldots,\widehat{s_i\circ f_{q,j}},\ldots, 
s_{p+1}\circ f_{q,j}).
\end{align*}
We also have
\begin{align*}
\delta(\dd{\phi})(s)
&= \delta\left( (-1)^{p+1}\sum_{j=0}^{q+1} (-1)^j {\phi}(
s_0\circ f_{q,j}, \ldots, s_{p+1}\circ f_{q,j})\right);\\
&=(-1)^{p+1} \sum_{j=0}^{q+1} \sum_{i=0}^{p+1} (-1)^{i+j} 
{\phi}(s_0\circ f_{q,j}, \ldots,\widehat{s_i\circ f_{q,j}},\ldots, s_{p+1}\circ f_{q,j}).
\end{align*}
Thus, it follows that $\dd\delta+\delta\dd=0$, so
$D^{\bullet,\bullet}$ is indeed a double complex.
\end{proof}

If $f: X \to Y$ is locally split, and if $\fU$ is the corresponding
open covering of $X$ defined in Section~\ref{subsec:singular}, the double
complex $D^{\bullet,\bullet}$ induces by restriction a double complex
$D_\fU^{\bullet,\bullet}$, where $D_\fU^{p,q}=C^q_\fU(W^p_f(X))$ when
$p\geq 0$ and $q \geq 0$ and $D_\fU^{\bullet,\bullet}=0$ otherwise.

The initial terms of 
the two spectral sequences associated with $D_\fU^{\bullet,\bullet}$ 
(cf. Section~\ref{subsec:spectral}) are as follows.
The first terms of the spectral sequence $'E_*^{i,j}$ are 
${'E}_1 = H_{\delta}(D_\fU^{\bullet,\bullet}(X)), 
{'E}_2 = H_{\dd}H_{\delta}(D_\fU^{\bullet,\bullet}(X)).
$
By the exactness of the sequence in Proposition~\ref{prop:MV}, 
we have that the spectral sequence 
${'E}_*^{i,j}$ degenerates at the $'E_2$ term as shown below.

\[
{'E}_1 =
\begin{array}{|cccccccccccccc}
& \vdots  && \vdots&& \vdots&& \vdots  && \vdots  && \\

&\Big\uparrow\vcenter{\rlap{$d$}} & &
\Big\uparrow\vcenter{\rlap{$0$}} & &
\Big\uparrow\vcenter{\rlap{$0$}} & &
\Big\uparrow\vcenter{\rlap{$0$}} & &
\Big\uparrow\vcenter{\rlap{$0$}}& & \\

& C_\fU^3(Y) &
 & 0 &
 & 0 &
 & 0 &
 & 0 & \cdots &
  
\\

&\Big\uparrow\vcenter{\rlap{$d$}} & &
\Big\uparrow\vcenter{\rlap{$0$}} & &
\Big\uparrow\vcenter{\rlap{$0$}} & &
\Big\uparrow\vcenter{\rlap{$0$}} & &
\Big\uparrow\vcenter{\rlap{$0$}} & & \\

& C_\fU^2(Y) &
 & 0 &
 & 0 &
 & 0 &
 & 0 & \cdots &
 
\\

&\Big\uparrow\vcenter{\rlap{$d$}} & &
\Big\uparrow\vcenter{\rlap{$0$}} & &
\Big\uparrow\vcenter{\rlap{$0$}} & &
\Big\uparrow\vcenter{\rlap{$0$}} & &
\Big\uparrow\vcenter{\rlap{$0$}} & & \\

& C_\fU^1(Y) &
 & 0 &
 & 0 &
 & 0 &
 & 0 & \cdots &
 
\\

&\Big\uparrow\vcenter{\rlap{$d$}} & &
\Big\uparrow\vcenter{\rlap{$0$}} & &
\Big\uparrow\vcenter{\rlap{$0$}} & &
\Big\uparrow\vcenter{\rlap{$0$}} & &
\Big\uparrow\vcenter{\rlap{$0$}} & & \\

&C_\fU^0(Y) &
 & 0 &
 & 0 &
 & 0 &
 & 0 & \cdots &
\\ \\
\hline
\end{array}
\]

and, by Proposition~\ref{prop:hatcher}, 
\[
'E_2 = 
\begin{array}{|cccccccc}
 & \vdots  &\vdots & \vdots  & \vdots & \vdots  & \vdots\\
& &
 & &
 & &
 &  \\

 & {H}^3(Y) & 0
 & 0 & 0
 & 0 & 0& \cdots 
\\

 &
 & &
 & &
 & & \\

  & {H}^2(Y) & 0
 & 0 & 0
 & 0 & 0 &\cdots
 
\\

 &
 & &
 & &
 & & \\

  & {H}^1(Y) & 0
 & 0 & 0
 & 0 & 0 &\cdots 
\\

 &
 & &
 & &
 & & \\

 & {H}^0(Y) & 0
 & 0 & 0
 & 0 & 0 &\cdots 
\\

& 
 & &
 & &
 & & \\
\hline
\end{array}
\]

The degeneration of this sequence at ${'E}_2$ shows that
\[
H^*(\Tot^{\bullet}(D_\fU^{\bullet,\bullet}(X))) \cong {H}^*(Y).
\]

The initial term ${''E}_1$ of the second spectral sequence is given by,
\[
{''E}_1 =
\begin{array}{|ccccccc}

 & \vdots  && \vdots  && \vdots  & \\
 &
 & &
 & &
 &  \\

 & {H}^3(W_f^0(X)) &
\stackrel{\delta}{\longrightarrow} & 
{H}^3(W_f^1(X)) &
\stackrel{\delta}{\longrightarrow} & 
{H}^3(W_f^2(X)) &
 \longrightarrow 
\\

 &
 & &
 & &
 &  \\

 & {H}^2(W_f^0(X)) &
\stackrel{\delta}{\longrightarrow} & 
{H}^2(W_f^1(X)) &
\stackrel{\delta}{\longrightarrow} & 
{H}^2(W_f^2(X)) &
 \longrightarrow 
\\

 &
 & &
 & &
 &  \\

 & {H}^1(W_f^0(X)) &
\stackrel{\delta}{\longrightarrow} & 
{H}^1(W_f^1(X)) &
\stackrel{\delta}{\longrightarrow} & 
{H}^1(W_f^2(X)) &
 \longrightarrow 
\\

 &
 & &
 & &
 &  \\

 & {H}^0(W_f^0(X)) &
\stackrel{\delta}{\longrightarrow} &
{H}^0(W_f^1(X)) &
\stackrel{\delta}{\longrightarrow} & 
{H}^0(W_f^2(X)) &
 \longrightarrow \\\\
\hline 
\end{array}
\]

Since this spectral sequence also converges to 
$H^*(\Tot^{\bullet}(D_\fU^{\bullet,\bullet})(X)),$ we have the
following proposition.

\begin{proposition}
\label{prop:intermediate}
$$
\displaylines{
H^*(\Tot^{\bullet}(D_\fU^{\bullet,\bullet})(X))  \cong {H}^*(Y).
}
$$ 
\end{proposition}

Proposition \ref{prop:intermediate} now implies,
\begin{theorem}
\label{the:descent}
For any continuous  semi-algebraic surjection $f:X \to Y$,
where $X$ and $Y$ are open semi-algebraic subsets of $\R^n$ and $\R^m$ 
respectively (or, more generally, for any locally split continuous
surjection $f$),
the spectral sequence associated to the double complex
$D^{\bullet,\bullet}(X)$ with 
$E_1 = H_d(D^{\bullet,\bullet}(X))$ converges  to 
$H^*(C^\bullet(Y)) \cong {H}^*(Y).$
In particular,
\begin{enumerate}
\item
$\displaystyle{
E_1^{i,j} = {H}^j(W_f^i(X)),
}
$ and
\item
$\displaystyle{
E_\infty \cong  H^*(\Tot^{\bullet}(D^{\bullet,\bullet}(X))) \cong {H}^*(Y).
}
$
\end{enumerate}
\end{theorem}

\begin{proof}
\hide{
The previous reasoning shows the result is true for the double
complex $D_\fU^{\bullet,\bullet}$. But by
Proposition~\ref{prop:hatcher}, it is clear that the total complexes
corresponding to the double complexes
$D_\fU^{\bullet,\bullet}$ and $D^{\bullet,\bullet}$ have isomorphic
cohomology. Thus, we have $E_\infty \cong
H^*(\Tot^{\bullet}(D_\fU^{\bullet,\bullet}(X))) \cong {H}^*(Y)$ (the
condition on $E_1$ follows immediately from the definition). 
}
By Proposition \ref{prop:hatcher}, we have that the 
component-wise homomorphisms, ${\iota}^{\bullet}$, induces a homomorphism
of double complexes,
\[
\iota^{\bullet,\bullet}: 
D^{\bullet,\bullet} \rightarrow D_{\fU}^{\bullet,\bullet},
\]
which in turn induces an isomorphism between the $E_1$ terms of the
corresponding spectral sequences.
Hence, by Theorem~\ref{the:spectral} we have that,
$H^*(\Tot^{\bullet}(D_{\fU}^{\bullet,\bullet})) \cong
H^*(\Tot^{\bullet}(D^{\bullet,\bullet})).$
The Theorem now follows from Proposition \ref{prop:intermediate}.
\end{proof}

\subsection{Truncation of the double complex}
\label{subsec:trunc}
If we  denote by $D_q^{\bullet,\bullet}(X)$  the truncated complex
defined by,
$$
\begin{array}{cccc}
D_q^{i,j}(X) &=& D^{i,j}(X), &  \;\;
\hbox{if} \quad 0 \leq i+j \leq q+1, \\
                    &=& 0,             & \;\;\mbox{otherwise},
\end{array}
$$
then it is clear that, 
\begin{equation}
\label{eqn:tot}
{H}^i(Y) \cong H^i(\Tot^{\bullet}(D_q^{\bullet,\bullet}(X))), 
\;\;\hbox{for}\quad  0 \leq i \leq q. 
\end{equation}

Now suppose that $X \subset \R^{k+m}$ is a compact semi-algebraic
set defined by the inequalities,
$P_1 \geq 0,\ldots,P_\ell \geq 0$.
Let $\pi$ denote the projection map,
$\pi : \R^{k+m}  \rightarrow \R^m$. 
Let $\epsilon > 0$ 
and let $\tilde{X} \subset \R^{k+m}$ be the set defined by $P_1 +
\epsilon > 0, \ldots, P_\ell + \epsilon > 0.$
\begin{proposition}
\label{prop:open}
\begin{enumerate}
\item For $\epsilon >0$ sufficiently small, we have 
$$
\displaylines{
{H}^*(W^p_\pi(\tilde{X})) \cong H^*(W^p_\pi(X)), \hbox{ for all } p
\geq 0,\cr 
\hbox{ and }{H}^*(\pi(\tilde{X})) \cong H^*(\pi(X)).
}
$$
\item
The map,
$\pi|_{\tilde{X}}$ is 
a locally split semi-algebraic surjection onto its image.
\end{enumerate}
\end{proposition}

\begin{proof}
  When $\epsilon >0$ is small, the sets $X$ and $\tilde{X}$ are homotopy
  equivalent and so are the sets $\pi(X)$ and $\pi(\tilde{X})$ and the
  fibered products $W^p_\pi(\tilde{X})$ and $W^p_\pi(X)$ for all $p
  \geq 0$ (see~\cite{B99}).  The first part of the proposition follows
  from the homotopy invariance property of singular cohomology groups.
  The second part of the proposition is clear once we note that
  $\tilde{X}$ is an open subset of $\R^{k+m}$: projections of open
  sets always admit local continuous sections.
\end{proof}

We can combine Theorem~\ref{the:descent} and
Proposition~\ref{prop:open} to construct, from the projection of a
{\em compact} basic semi-algebraic set, a double complex giving rise
to a cohomological descent spectral sequence.
\begin{corollary}
\label{cor:descent2}
Let $X \subset \R^{k+m}$ be a compact semi-algebraic set defined
by $P_1 \geq 0, \ldots,P_\ell \geq 0$ and $\pi:R^{k+m} \rightarrow \R^m$
the projection onto the last $m$ co-ordinates.
The spectral sequence associated to the double complex
$D^{\bullet,\bullet}(X)$ with 
$E_1 = H_d(D^{\bullet,\bullet}(X))$ converges  to 
$H^*(C^\bullet(\pi(X))) \cong {H}^*(\pi(X)).$
In particular,
\begin{enumerate}
\item
$\displaystyle{
E_1^{i,j} = {H}^j(W_f^i(X)),
}
$
and
\item
$
\displaystyle{
E_\infty \cong  H^*(\Tot^{\bullet}(D^{\bullet,\bullet}(X))) \cong 
{H}^*(\pi(X)).
}
$
\end{enumerate}
\end{corollary}

\begin{remark}
  Note that it is
  not obvious how to prove directly an exact sequence at the level of singular
  (or even simplicial) cochains for the projection of a compact set, 
  as we do in Proposition~\ref{prop:MV} in the locally-split setting. 
  One difficulty is the fact  that semi-algebraic maps are not, in general, 
  triangulable.
\end{remark}

Now let $X$ be a compact semi-algebraic set defined
by a constant number of quadratic inequalities and $f$ a projection map. 
We cannot hope to compute even the truncated complex
$D_q^{\bullet,\bullet}(X)$ since these are defined in terms of
singular chain complexes 
which are infinite-dimensional.
We overcome this problem by 
computing  another double complex
$\D_q^{\bullet,\bullet}(X)$,
such that there exists a homomorphism of double complexes,
$\displaystyle{
\psi: \D_q^{\bullet,\bullet}(X) \longrightarrow 
D_q^{\bullet,\bullet}(X),
}
$
which induces an isomorphism between the $'E_1$ terms of the
spectral sequences associated to 
the double complexes $D_q^{\bullet,\bullet}(X)$ and 
$\D_q^{\bullet,\bullet}(X).$ 
This implies, by virtue of Theorem~\ref{the:spectral}, 
that the cohomology groups of the associated total complexes are isomorphic,
that is,
$$
\displaylines{
H^*(\Tot^{\bullet}(D_q^{\bullet,\bullet}(X)))
\cong 
H^*(\Tot^{\bullet}(\D_q^{\bullet,\bullet}(X))).
}
$$

The construction of the double complex 
$\D_q^{\bullet,\bullet}(X)$ is 
described in Section~\ref{sec:main}.

\section{Algorithmic Preliminaries}
\label{sec:algo}

We now recall an algorithm described in \cite{Basu05a}, where the following
theorem is proved.

\begin{theorem}
\label{the:main2}
There exists an algorithm, which takes as input 
a family of polynomials
$\{P_1,\ldots,P_s\}\subset  \R[X_1\ldots,X_k],$ with
${\rm deg}(P_i) \leq 2,$
and a number $\ell \leq k$,
and outputs a complex ${\mathcal D}^{\bullet,\bullet}_\ell.$ 
The complex $\Tot^{\bullet}({\mathcal D}^{\bullet,\bullet}_\ell)$ is
quasi-isomorphic to $\C^\ell_{\bullet}(S)$, the truncated singular
chain complex of $S$,
where
$$S = \bigcap_{P \in {\mathcal P}}
               \{x \in \R^{k}\; \mid \; P(x) \leq 0 \}.
$$

Moreover, given a subset ${\mathcal P}' \subset {\mathcal P}$, 
with 
$$
S' = \bigcap_{P \in {\mathcal P}'}
               \{x \in \R^{k}\; \mid \; P(x) \leq 0 \}.
$$
the algorithm outputs both complexes ${\mathcal D}^{\bullet,\bullet}_\ell$ and
${\mathcal D}'^{\bullet,\bullet}_\ell$ (corresponding to the sets
$S$ and $S'$ respectively) along with the matrices defining 
a homomorphism $\Phi_{{\mathcal P},{\mathcal P}'},$ 
such that 
$\Phi_{{\mathcal P},{\mathcal P}'}^*: H^*(\Tot^\bullet({\mathcal D}^{\bullet,\bullet}_\ell)) \cong H^*(S) \rightarrow  H^*(S') \cong
H^*(\Tot^\bullet({\mathcal D'}^{\bullet,\bullet}_\ell))
$
is the homomorphism induced by the inclusion $i: S \hookrightarrow S'.$
The complexity of the algorithm is 
$ 
\sum_{i=0}^{\ell+2} %
\binom{s}{i}
k^{2^{O(\min(\ell,s))}}.
$
\end{theorem}

For completeness, we formally state
the input and output of the algorithm mentioned in
Theorem \ref{the:main2}.

We first introduce some notations
which will be used to describe the input and output of the algorithm.
Let
${\mathcal Q} = \{Q_1,\ldots,Q_s \} \subset  \R[X_1,\ldots,X_k]$ be a family of
polynomials
with $\deg(Q_i) \leq 2, 1 \leq i \leq s$.
For each subset $J \subset \{1,\ldots,s\}$,
we let $S_J$ denote the semi-algebraic
set defined by %
$\{Q_j \geq 0 \mid j \in J\}$.
Notice that for each pair $I \subset J \subset \{1,\ldots,s\},$
we have an inclusion $S_J \subset S_I$.
\begin{algorithm}[Build Complex]
\label{alg:complex}
\item [{\bf Input:}] 
A family of polynomials
${\mathcal Q} = \{Q_1,\ldots,Q_s \} \subset  \R[X_1,\ldots,X_k]$
with $\deg(Q_i) \leq 2$, for  $1 \leq i \leq s$.
\item [{\bf Output:}]
\begin{enumerate}
\item[]
\item
For each subset $J \subset \{1,\ldots,s\}$,
a description of a complex $F^{\bullet}_J$, 
consisting of a basis for each term of the
complex and matrices (in this basis) for the differentials, and
\item
for each pair $I \subset J \subset \{1,\ldots, s\}$,
a homomorphism, $\phi_{I,J}: F^{\bullet}_I \longrightarrow F^{\bullet}_J.$
\end{enumerate}
\noindent The complexes, $F_J^{\bullet}$ and the homomorphisms $\phi_{I,J}$
satisfy the following.
\begin{enumerate}
\item
For each $J \subset \{1,\ldots,s\}$,
\begin{equation}
\label{eqn:iso}
H^*(F^{\bullet}_J) \cong H^*(S_J).
\end{equation}

\item
For each pair $I \subset J \subset \{1,\ldots,s\},$
the following diagram commutes.
\begin{diagram}
H^*(F^{\bullet}_I)& \rTo^{\left(\phi_{I,J}\right)^*}&
H^*(F^{\bullet}_J) \\
\uTo^{\cong}&&\uTo^{\cong}\\
H^*(S_I) &\rTo^{r^*} &H^*(S_J)
\end{diagram}
Here, $\left(\phi_{I,J}\right)^*$ is the homomorphism induced by
$\phi_{I,J},$ the vertical homomorphisms are the isomorphisms from
\eqref{eqn:iso}, and $r^*$ is the homomorphism induced by 
restriction.
\end{enumerate}
\end{algorithm}

\noindent
{\rm Complexity:}
The complexity of the algorithm is $k^{2^{O(s)}}.$
\eop

For the purposes of this paper, we need to slightly modify 
Algorithm~\ref{alg:complex} 
in order 
to be able to handle permutations of the co-ordinates.
More precisely, suppose that
$\sigma \in \mathfrak{S}_k$ is a given permutation of the co-ordinates,
and for any $I \subset \{1,\ldots,s\},$
let $S_{I,\sigma} = \{ (x_{\sigma(1)},\ldots, x_{\sigma(k)})\; \mid\;
(x_1,\ldots,x_k) \in S_I \}.$
Let $F_{I,\sigma}^\bullet$ denote  the complex computed by the
algorithm corresponding to the set $S_{I,\sigma}$.
It is easy to modify 
Algorithm~\ref{alg:complex} slightly 
without changing the complexity estimate,
such that for any fixed $\sigma$,
the algorithm outputs, complexes $F_I^\bullet,F_{I,\sigma}^\bullet$
as well as the matrices corresponding to the 
induced isomorphisms,
$\phi^{\bullet}_\sigma: F_I^\bullet \rightarrow F_{I,\sigma}^\bullet.$
We assume this implicitly in the description of Algorithm
\ref{alg:main} in the next section.

\section{Algorithm for projections}
\label{sec:main}
Let $S \subset \R^{k+m}$ be a basic semi-algebraic set defined by 
\[
P_1 \geq  0, \ldots, P_\ell \geq 0,
P_i \in \R[X_1,\ldots,X_k,Y_1,\ldots,Y_m],
\]
with $\deg(P_i) \leq 2, \; 1 \leq i \leq \ell$. 
Let $\pi:\R^{k+m} \rightarrow \R^m$ be the projection onto the last
$m$ coordinates.

The algorithm will compute a double complex, 
$\D^{\bullet,\bullet}_q(S)$, such that
$\Tot^{\bullet}(\D^{\bullet,\bullet}_q(S))$ 
is quasi-isomorphic to the complex
$\Tot^{\bullet}(D^{\bullet,\bullet}_q(S))$.
The double complex, $\D_q^{\bullet,\bullet}(S)$ is defined
as follows.

We introduce $k(q+2)$ variables, which we denote by
$X_{i,j}$, $1 \leq i \leq k, 0 \leq j \leq q+1$. 
For each $j, 0 \leq j \leq q+1$,
we denote by,
$P_{i,j}$ the polynomial 
$$P_i(X_{1,j},\ldots,X_{k,j},Y_1,\ldots,Y_m)$$
(substituting $X_{1,j},\ldots,X_{k,j}$ in place of
$X_1,\ldots,X_k$ in the polynomial $P_i$).
We consider each $P_{i,j}$ to be an element of
$\R[X_{1,0},\ldots,X_{k,q+1},Y_1,\ldots,Y_m]$.
For each $p, 0 \leq p \leq q+1$, we denote by 
$S_p \subset \R^{k(q+2)+m}$ the semi-algebraic set defined by,
\[
P_{1,0}\geq 0,\ldots,P_{\ell,0} \geq 0, \ldots,
P_{1,p} \geq 0,\ldots,P_{\ell,p} \geq 0.
\]

Note that, for each $p, 0 < p \leq q+1$, 
and each $j, 0 \leq j \leq  p$ we have a natural 
map,
$\pi_{p,j}: S_{p} \rightarrow S_{p-1}$ given by,
$$
\pi_{p,j}(\bar{x}_0,\ldots,\bar{x}_p,\ldots,\bar{x}_{q+1},\bar{y})
= (\bar{x}_0,\ldots,\bar{x}_p,\ldots,\bar{x}_j,\ldots,\bar{x}_{q+1},\bar{y}).
$$
Note that in the definition above, each $\bar{x}_i \in \R^k$ and
$\pi_{p,j}$ exchanges the coordinates $\bar{x}_j$ and
$\bar{x}_p$.

We are now in a position to define $\D_q^{\bullet,\bullet}$.
We follow the notations introduced in Section~\ref{sec:algo}.
Let ${\mathcal Q}= \{Q_1,\ldots,Q_{\ell(q+2)}\}= 
\{P_{1,0},\ldots,P_{\ell,q+1}\}.$
For $0 \leq j \leq q+1$, we let $L_j = \{1,\ldots,(j+1)\ell\} \subset
\{1,\ldots,(q+2)\ell\}.$

$$
\begin{array}{cccc}
\D_q^{i,j}(X) &=& F^{j}_{L_i}, &  \;\;
0 \leq i + j\leq q+1, \\
                    &=& 0,             & \;\;\mbox{otherwise},
\end{array}
$$

The vertical homomorphisms, $d,$ in the complex $\D_q^{\bullet,\bullet}$
are those induced from the complexes $F_{L_i}^{\bullet}$ or zero.
The horizontal homomorphisms,
$\delta^j: F_{L_{i}}^j \longrightarrow F_{L_{i+1}}^j$ are defined as follows.

For each $h, 0 \leq h \leq i+1$, Algorithm~\ref{alg:complex}
produces a homomorphism,
$\phi_{i+1,h}: F_{L_{i}}^j \longrightarrow F_{L_{i+1}}^j$, corresponding to
the map $\pi_{i+1,h}$
(see remark after Algorithm~\ref{alg:complex}).
The homomorphism $\delta$ is then defined
by, 
$
\displaystyle{
\delta = \sum_{h=0}^{i+1} (-1)^h \phi_{i+1,h}.
}
$
We have the following proposition.
\begin{proposition}
\label{prop:correctness}
The complex
$\Tot^{\bullet}(\D^{\bullet,\bullet}_q(S))$ 
is quasi-isomorphic to the complex
$\Tot^{\bullet}(D^{\bullet,\bullet}_q(S))$.
\end{proposition}
\begin{proof}
It follows immediately from 
Theorem \ref{the:main2} that the columns of the complexes
$\D^{\bullet,\bullet}_q(S)$  and
${D}^{\bullet,\bullet}_q(S)$ are quasi-isomorphic.
Moreover, it is easy to see that the quasi-isomorphisms
induce an isomorphism between the $''E_1$ term of their
associated spectral sequences. Now by Theorem \ref{the:spectral} this
implies that 
$\Tot^{\bullet}(\D^{\bullet,\bullet}_q(S))$ 
is quasi-isomorphic to the complex
$\Tot^{\bullet}(D^{\bullet,\bullet}_q(S))$.
\end{proof}

\begin{algorithm} [Computing the first $q$  Betti Numbers]
\label{alg:main}
\item [{\bf Input:}] 
A $S \subset \R^{k+m}$ be a basic semi-algebraic set defined by 
\[
P_1 \geq  0, \ldots, P_\ell \geq 0,
\] 
with $P_i \in \R[X_1,\ldots,X_k,Y_1,\ldots,Y_m],$
$\deg(P_i) \leq 2, \; 1 \leq i \leq \ell$. 

\item [{\bf Output:}] 
$b_{0}(\pi(S)),\ldots, b_q(\pi(S))$,
where
$\pi:\R^{k+m} \rightarrow \R^m$ be the projection onto the last
$m$ coordinates. 

\item [{\bf Procedure:}]
\item[Step 1:]
Using Algorithm~\ref{alg:complex} compute the 
truncated complex $\D_q^{\bullet,\bullet}(S)$.

\item[Step 2:]
Compute using linear algebra, the dimensions of 
$\displaystyle{
H^i(\Tot^{\bullet}(\D_q^{\bullet,\bullet})), 
\;\;0 \leq i \leq q. 
}
$

\item[Step 3:]
For each $i, \;\;0 \leq i \leq q,$
output, 
$b_i(\pi(S)) = 
\dim(H^i(\Tot^{\bullet}(\D_q^{\bullet,\bullet}))).$
\end{algorithm}

\noindent
{\rm Complexity Analysis:}
The calls to Algorithm~\ref{alg:complex} has input consisting of
$(q+1)\ell$ polynomials in $qk+m$ variables. Using the complexity
bound of Algorithm~\ref{alg:complex} we see that the complexity
of Algorithm~\ref{alg:main} is bounded by $(k + m)^{2^{O(q\ell)}}$
\eop

\noindent
{\rm Proof of Correctness:}
The correctness of the algorithm is a consequence of 
Proposition \ref{prop:correctness}
and Theorem~\ref{the:spectral}.
\eop

\section{Conclusion and Open Problems}
For any fixed $q$ and $\ell$, we have proved a polynomial bound on the
sum of the first $q$ Betti numbers of the projection of a bounded, basic
closed semi-algebraic set defined by $\ell$ quadratic inequalities. 
We have also described a polynomial time algorithm to compute the 
first $q$ Betti numbers of the image of such a projection.

Since it is not known whether quantifier elimination can be performed
efficiently for sets defined by a fixed number of quadratic
inequalities, many questions are left open.

Our bounds become progressively worse as $q$ increases, 
becoming exponential in the dimension as $q$
approaches $k$. However, we do not have any examples (of projections of
semi-algebraic sets defined by quadratic inequalities) where  the
higher Betti numbers behave exponentially in the dimension. This leaves open
the problem of either constructing such examples, or removing the dependence
on $q$ from our bounds.

Another interesting open problem  is to improve the complexity of 
Algorithm~\ref{alg:main}, 
from $(k+m)^{2^{O(q\ell)}}$ to $(k+m)^{O(q\ell)}$. Note
that this would imply an algorithm with complexity
$k^{O(q\ell)}$ for computing the first $q$ Betti numbers of a 
semi-algebraic set defined by $\ell$ quadratic inequalities in $\R^k$. 
The best known algorithm for computing all the
Betti numbers of such sets has complexity $k^{2^{O(\ell)}}$~\cite{Basu05a}. 
The only topological invariants of such sets
that we currently know how to compute in time $k^{O(\ell)}$ are
testing for emptiness  \cite{Barvinok93,GP} and the Euler-Poincar\'e
characteristic~\cite{Basu05c}.


\begin{thebibliography}{50}

\bibitem{Agrachev}
{\sc A. A.\ Agrachev},
\newblock {\em Topology of quadratic maps and Hessians of smooth maps},
\newblock {Algebra, Topology, Geometry, Vol 26 (Russian),85-124, 162,
Itogi Nauki i Tekhniki, Akad. Nauk SSSR, Vsesoyuz. Inst. Nauchn.i 
Tekhn. Inform., Moscow, 1988.
Translated in J. Soviet Mathematics. 49 (1990), no. 3, 990-1013.}

\bibitem{Barvinok93}
{\sc A. I.\ Barvinok},
\newblock {\em Feasibility Testing for Systems of Real Quadratic Equations},
\newblock {Discrete and Computational Geometry}, 10:1-13 (1993).

\bibitem{Barvinok97}
{\sc A. I.\ Barvinok}
\newblock {\em On the Betti numbers of semialgebraic sets defined by few
quadratic inequalities}
\newblock {\em Math. Zeit.}, 225:231-244, 1997.


\bibitem{B99}
{\sc S.\ Basu},
\newblock {\em On Bounding the Betti Numbers and Computing the Euler
Characteristics of Semi-algebraic Sets},
\newblock {Discrete and Computational Geometry}, 22 1-18 (1999).

\bibitem{BPR95}
{\sc S.\ Basu, R.\ Pollack, M.-F. \ Roy},
 \newblock {\em On the Combinatorial and Algebraic
Complexity of Quantifier Elimination},
\newblock Journal of the ACM, 43  1002--1045, (1996).


\bibitem{B03}
{\sc S.\ Basu},
\newblock{\em On different bounds on different Betti numbers},
\newblock {Discrete and Computational Geometry}, 30:1, 65-85, 2003.


\bibitem{Basu05a}
{\sc S.\ Basu},
\newblock{\em Polynomial time algorithm for computing the top Betti numbers
of Semi-algebraic sets defined by quadratic inequalities},
\newblock {Proceedings of Symposium on the Theory of Computing}, 2005.
Available at 
{\tt www.math.gatech.edu/$\sim$saugata/quadratic.ps}.

\bibitem{Basu05b}
{\sc S.\ Basu},
\newblock{\em Single exponential time algorithm for computing the first 
few Betti numbers of semi-algebraic sets},
\newblock preprint.
Available at 
{\tt www.math.gatech.edu/$\sim$saugata/bettifew.ps}

\bibitem{Basu05c}
{\sc S.\ Basu},
\newblock{\em 
Efficient algorithm for computing the Euler-Poincar\'e characteristic
of semi-algebraic sets defined by few quadratic inequalities},
\newblock to appear in Computational Complexity. 
Available at 
{\tt www.math.gatech.edu/$\sim$saugata/eulerquad.pdf}


\bibitem{BPR05}
{\sc S.\ Basu, R.\ Pollack, M.-F. \ Roy},
\newblock{\em Computing  the first Betti number and the connected components
of semi-algebraic sets},
\newblock {Proceedings of Symposium on the Theory of Computing}, 2005 
(to appear).
Available at 
{\tt www.math.gatech.edu/$\sim$saugata/bettione.ps}.

\bibitem{BPR03}
{\sc S.\ Basu, R.\ Pollack, M.-F. \ Roy},
\newblock {\em Algorithms in Real Algebraic Geometry},
Springer-Verlag, 2003.

\bibitem{BCR}
{\sc J.\ Bochnak, M.\ Coste, M.-F.\ Roy},
\newblock {\em G\'eom\'etrie
alg\'ebrique r\'eelle,}
Springer-Verlag (1987).
\newblock {\em Real algebraic Geometry},
Springer-Verlag (1998).


\bibitem{BT}
{\sc R.\ Bott, L. W.\ Tu},
\newblock {\em Differential Forms in Algebraic Topology},
Springer-Verlag (1982).

\bibitem{BC}
{\sc P.\ Burgisser, F.\ Cucker},
\newblock {\em Counting Complexity Classes for Numeric Computations II: 
Algebraic and  Semi-algebraic Sets},
\newblock preprint.

\bibitem{Bredon}
{\sc G. E. Bredon},
\newblock {\em Sheaf Theory,}
Springer-Verlag (1996).

\bibitem{Col}
{\sc G. Collins},
\newblock {\em Quantifier elimination for real closed fields by
cylindric algebraic decomposition},
\newblock In Second GI Conference on Automata Theory and
Formal Languages. Lecture Notes in Computer Science, vol. 33, pp. 134-183, 
Springer-Verlag, Berlin (1975).

\bibitem{deligne}
{\sc P. Deligne},
\newblock {\em Th\'eorie de Hodge~III},
\newblock { Publ. Math. IHES} 44:5--77, 1974.


\bibitem{dugger}
{\sc D.~Dugger} and {\sc D.~Isaksen},
\newblock {\em Topological hypercovers and $A^1$-realizations},
\newblock Math. Z 246 (2004), 667--689.

\bibitem{GVZ}
{\sc A.\ Gabrielov, N.\ Vorobjov, T.\ Zell},
\newblock {\em Betti Numbers of Semi-algebraic and Sub-Pfaffian Sets},
\newblock {J. London Math. Soc. (2) 69 (2004) 27-43.}

\bibitem{Gabrielov}
{\sc A.\ Gabrielov},
\newblock{\em Counter-examples to quantifier elimination for fewnomial and 
exponential expressions, preprint.}
\newblock{Available at {\tt
http://www.math.purdue.edu/$\sim$agabriel/preprint.html}}.




\bibitem{GP}
{\sc D.\ Grigor'ev, D.V.\ Pasechnik},
\newblock {\em Polynomial time computing over quadratic maps I. 
Sampling in real algebraic sets},
\newblock Computational Complexity 14:20-52, (2005).


\bibitem{Hardt} {\sc R. M. Hardt},
\newblock{\em Semi-algebraic Local Triviality in Semi-algebraic Mappings},
\newblock Am. J. Math. { 102}, 291-302 (1980).

\bibitem{Hatcher} {\sc A. Hatcher},
\newblock{\em Algebraic Topology},
\newblock Cambridge University Press (2002).

\bibitem{houston}
{\sc K.~Houston},
\newblock {\em An introduction to the image computing spectral sequence},
\newblock In {\em Singularity theory (Liverpool, 1996)}, volume 263 of {\em
  London Math. Soc. Lecture Note Ser.}, pages 305--324. Cambridge
  Univ. Press, Cambridge, 1999.


\bibitem{Khovansky}
{\sc A. G.\ Khovansky}
\newblock{Fewnomials,}
\newblock{American Mathematical Society}, 1991.


\bibitem{Mcleary}
{\sc J. \ McCleary}
\newblock A User's Guide to Spectral Sequences, Second Edition 
\newblock Cambridge Studies in Advanced Mathematics, 2001.

\bibitem{Milnor}
{\sc J. \ Milnor},
\newblock {\em On the Betti numbers of real varieties},
\newblock Proc. AMS 15, 275-280 (1964).

\bibitem{murray}
{\sc M. Murray},
\newblock {\em Bundle gerbes},
\newblock {J. London Math. Soc.} {54}, 403-416 (1996).



\bibitem{O}
{\sc O.\ A.\ Ole\u{i}nik},
\newblock {\em Estimates of the {B}etti numbers
of real algebraic hypersurfaces},
\newblock { Mat. Sb. (N.S.)}, 28 (70): 635--640 (Russian) (1951).
%


\bibitem{OP}
{\sc O. A.\ Ole\u{i}nik, I. B.\ Petrovskii},
\newblock {\em On the topology of real algebraic surfaces},
\newblock Izv. Akad. Nauk SSSR 13, 389-402 (1949).
%

\bibitem{R92}
{\sc J.\ Renegar.}
\newblock {\em On the computational complexity and geometry of the
first order theory of the reals},
\newblock Journal  of Symbolic Computation, 13: 255--352 (1992).

\bibitem{stdonat}
{\sc B.~Saint-Donat},
\newblock {\em Techniques de descente cohomologique},
\newblock In {Th\'eorie des topos et cohomologie \'etale des sch\'emas.} 
Tome 2, Springer-Verlag, Berlin, 1972 (SGA 4), 
Lecture Notes in Mathematics, Vol. 270, p. 83--162.

%
%
%
%

\bibitem{Thom}
{\sc R.\ Thom},
\newblock {\em Sur l'homologie des vari\'et\'es alg\'ebriques
r\'eelles},
\newblock  {Differential and Combinatorial
Topology},  255--265.
Princeton University Press, Princeton (1965).

\bibitem{vassiliev}
{\sc V.~Vassiliev},
\newblock { Complements of discriminants of smooth maps: topology and
  applications}, volume~98 of {\em Translations of Mathematical Monographs}.
\newblock American Mathematical Society, Providence, RI, 1992.

%
%
%
%
%

\bibitem{Hausdorff}
{\sc T.\ Zell},
\newblock {\em Topology of definable Hausdorff limits},
\newblock Discrete Comput. Geom. 33, 423--443 (2005).

\end{thebibliography}
\end{document}